\input amstex
\documentstyle{amsp2}
\TagsOnRight
\magnification=1100
\define\g{{\goth g}}
\define\h{{\goth h}}

\define\k{{\goth k}}
\redefine\l{{\goth l}}

\define\m{{\goth m}}

\redefine\b{{\goth b}}

\define\a{{\goth a}}
\redefine\c{{\goth c}}
\define\={\overset\text{def}\to=}
\redefine\B{{\Cal B}}
\define\C{{\Bbb C}}
\define\R{{\Bbb R}}
\define\Z{{\Bbb Z}}

\redefine\D{{\Cal D}}

\TagsOnRight
\topmatter
\title
Compact homogeneous CR manifolds
\endtitle

\author
Dmitry V. Alekseevsky and Andrea F. Spiro\\
\phantom{aa}
\endauthor

\address
\phantom{ }\newline D. V. Alekseevsky\newline Department of
Mathematics
\newline Hull University
\newline gen. Cottingham Road, HU6 7RX Hull \newline
ENGLAND\newline \phantom{ }\newline \phantom{ }
\endaddress

\email d.v.alekseevsky\@maths.hull.ac.uk
\endemail

\address
\phantom{ }\newline
A. F. Spiro\newline
Dipartimento di Matematica e Fisica\newline
Universit\`a di Camerino\newline
62032 Camerino (Macerata)\newline
ITALY\newline
\phantom{ }\newline
\phantom{ }
\endaddress

\email
spiro\@campus.unicam.it
\endemail


\leftheadtext{D. V. Alekseevsky and A. F. Spiro}
\rightheadtext{Compact homogeneous
CR manifolds}

\abstract We classify all compact simply connected homogeneous CR
manifolds $M$ of codimension one and with non-degenerate Levi form up to
CR equivalence. The classification is based on our previous results and
on a description of the maximal connected compact group $G(M)$ of automorphisms
of $M$. We characterize also the standard homogeneous CR manifolds as
the homogeneous CR manifolds whose group $G(M)$ in not semisimple.
\endabstract

\subjclass Primary 32V05; Secondary 53D10  57S25 53C30
\endsubjclass
\keywords Homogeneous CR manifolds, Real Hypersurfaces, Contact Homogeneous
manifolds
\endkeywords

\endtopmatter
\document

\subhead 1. Introduction
\endsubhead
\bigskip
In a previous paper (\cite{AS}), we
classified all simply connected
compact homogeneous CR manifolds  $(M = G/L, \D, J)$
of a compact Lie group $G$ up to a $G$-equivariant isomorphism.
Here $(\D, J)$ is a $G$-invariant CR structure
on the homogeneous manifold $M = G/L$, where $\D$ is a
codimension one $G$-invariant distribution on $M$ and
$J$ is a complex structure on $\D$, which satisfies the integrability condition
(2.1) (see \S 2.1, below).  Note that by  the results in \cite{AHR}
and \cite{Sp},  any simply connected compact homogeneous
CR manifold admits a compact transitive group of automorphisms and hence
it can be represented in the above form. \par
In the present paper, we study when two such homogeneous CR manifolds
$(M = G/L, \D, J)$ and  $(M' = G'/L', \D', J')$ are CR equivalent, that is
when there is
a diffeomorphism $\phi: M \longrightarrow M'$ such that
$\phi_*(\D) \subset \D'$ and $\phi_*(Jv) = J' \phi_*(v)$ for all $v\in \D$.\par
This question is reduced to the description of maximally compact
connected group  of automorphisms of a homogeneous
CR manifold $(M = G/L,\D, J)$. We give such description using
a result by A. L. Onishchik about
the maximal compact groups  of holomorphic transformations of
flag manifolds.\par
We shortly recall  the main results in \cite{AS} about the
classification of compact simply connected homogeneous CR manifolds
$(M = G/L, \D, J)$ of a compact Lie group $G$. Such manifolds
are subdivided into three natural disjoint classes:
\roster
\item"a)" the {\it standard  homogeneous CR
manifolds\/}, that is  homogeneous $S^1$-bundles over a flag manifold $F$,
with the CR structure induced by an invariant complex structure on $F$; \par
\item"b)" the {\it Morimoto-Nagano spaces\/}, i.e.  the
sphere bundles $S(N)\subset TN$ of a compact
rank one symmetric space $N = G/H$, with the CR structure induced by the
natural complex structure of $TN = G^\C/H^\C$; \par
\item"c)" the manifolds which  admit a  non-trivial holomorphic fibration over
a flag manifold $(F,J_F)$
with typical fiber $S(S^k)$, where $k = 2, 3, 5, 7$ or $9$, respectively;
these manifolds are  $SU_n/T^1\cdot SU_{n-2}$,
$SU_p\times SU_q/T^1 \cdot U_{p-2}\cdot U_{q-2}$,
$SU_n/T^1\cdot SU_2\cdot SU_2\cdot SU_{n-4}$,
$SO_{10}/T^1\cdot SO_6$ and $E_6/T^1\cdot SO_8$.
\endroster
In this last case, the invariant CR structure
is determined by the invariant complex
structure  $J_F$ on $F$ and by
an invariant CR structure on the typical fiber,
which depends on one complex parameter.\par
\medskip
First of all, we prove that  a non-standard homogeneous CR
manifold (i.e. a manifold from class b) or c)) is never CR
equivalent to a standard CR manifold. Moreover if two non-standard
homogeneous CR manifolds $M = G/L$ and $M' = G'/L'$ are CR
equivalent, then either $M = \text{Spin}_7/\text{SU}_3$ and $M' =
\text{SO}_8/\text{SO}_6$ and they are both CR equivalent to a
sphere bundle $S(S^7) \subset T S^7$, or they are equivalent as
homogeneous manifolds, that is there exists an isomorphism
$\phi: G\to G'$ such that $\phi(L) = L'$. \par
 Moreover, as we
proved in  \cite{AS},  the CR structures of a non-standard
homogeneous CR manifold $M$, with a fixed underlying contact
distribution, are naturally parameterized by the points of the
unit disc $D\subset \C$. We show  here that two CR structures
corresponding  to $t, t' \in D$ are CR equivalent if and only if
$|t| = |t'|$. \par
\medskip
Now, let $M = G/L$ be a standard CR manifold and $\pi: M = G/L \to F = G/K$
 the associated holomorphic $S^1$-fibration over a flag manifold.\par
We prove that any
maximal connected compact Lie group $A$ of automorphisms
of $M$, which contains $G$, preserves
the holomorphic fibration $\pi: M \to F$. In particular $A$ acts on the
flag manifold $F$ as a group of holomorphic transformations.
 Conversely,
any maximal connected compact group $A$ of holomorphic transformations of the flag
manifold $F = G/K$, which contains $G$, acts on
$M = G/L$ as a maximal  compact, connected
semisimple group of CR transformations. \par
Therefore the construction of a maximal  compact semisimple group
of CR  transformations of $M = G/L$ reduces to the description
of the maximal compact group of holomorphic transformations of the flag manifold
$F$. This problem was solved  by Onishchik  in  \cite{On} (see Theorem 4.1).
 In particular, he discovered that there exist only  few
irreducible flag manifolds which admit two different transitive groups
of holomorphic transformations, namely $\C P^{2\ell-1}$ ($\ell>1$),
$Gr_2(\R^7)$ and $Com(\R^{2\ell +2})$ ($\ell >2$) (see Table 1).\par
Using Onishchik's result, we describe the maximal compact semisimple Lie group
$A \supset G$ of a given standard CR manifold $M = G/L$ and we prove the
following.
Let  $M = G/L$ and $M' = G'/L'$ be two standard  CR manifolds and
$M = A/B$ and $M' = A'/B'$ their representations as homogeneous
spaces  of the maximal compact semisimple Lie groups  $A\supset G$
and $A' \supset G'$.  Then $M$ is CR equivalent to $M'$ if and only if
the homogenous manifolds $A/B$ and $A'/B'$ are equivalent and the associated
flag manifolds $F = G/K = A/C$ and $F' = G'/K' = A'/C'$ are
$A$-equivariantly biholomorphic.\par
In particular, we obtain that if the flag manifold $F = G/K$, associated with a
standard CR manifold $M = G/L$, has no factor isomorphic to
$\C P^{2\ell-1}$ ($\ell>1$),
$Gr_2(\R^7)$ or $Com(\R^{2\ell +2})$ ($\ell >2$), then
$G$ is a maximal connected compact
semisimple group of CR transformations of $M$.\par
\medskip
We also prove that a maximal connected compact group $A$
of CR transformations
of a compact homogeneous CR manifold $M$ is semisimple
if and only if $M$ is non-standard. For a standard CR manifold $M$,
the group $A$ has a 1-dimensional center, which acts trivially
on the associated flag manifold. This gives a group characterization
of standard CR manifolds.\par
\medskip
As a final remark, we would like to stress the fact that, if the Levi form
of a compact  CR manifold $M$ is positive definite, then the full group
of CR transformations is non-compact if and only if $M  = S^{2n+1}$ (\cite{We}).
We do not know
which of the classified compact homogeneous CR manifolds (with indefinite
Levi form) has non-compact full group of automorphisms. \par
Note that in \cite{Ya}, Yamaguchi classified  homogeneous Levi non-degenerate
CR manifolds with sufficiently large group of automorphisms. These manifolds
are either compact quadrics
or quadrics with some  points deleted and all
of them have a non-compact group of automorphism.\par
Also,
Fu, Isaev and Krantz (\cite{FIK}) and Zwonek (\cite{Zw})
found examples of non-homogeneous compact CR manifolds
of codimension one with non-compact group of automorphisms.
In these examples the Levi form is indefinite and degenerates at
some points.
\bigskip
\bigskip
\subhead 2. Preliminaries
\endsubhead
\bigskip
\subsubhead 2.1 First definitions
\endsubsubhead
\medskip
A  CR structure on a manifold $M$ is a
pair $(\Cal D , J)$, where $\Cal D \subset TM$ is a distribution on
$M$ with a complex structure $J$,
that is a field of endomorphisms
 $J\in \operatorname{End}\Cal D$, with $J^2 = -1$.
A CR structure $(\Cal D, J)$ is called {\it integrable\/} if $J$
satisfies the following integrability condition:
$$J([JX, Y] + [X, JY])\in \D\ ,$$
$$[JX, JY] - [X, Y] - J([JX, Y] + [X, JY])  = 0\ \tag2.1$$
for any pair of vector fields $X$, $Y$ in $\D$.\par
Geometrically this means that the eigendistributions
$\D^{10} \subset T^\C M$ and $\D^{01} \subset T^\C M$ of $J$, given
by the $J$-eigenspaces in $\D^\C$ corresponding to the eigenvalues $i$ and $-i$,
are involutive, i.e.  the space of their local sections is closed
under Lie brackets.\par
The {\it codimension} of a CR structure $(\Cal D, J)$
is defined as the codimension of the distribution $\Cal D$.
 An integrable
 codimension one CR structure $(\Cal D, J)$ is often called
{\it CR structure  of hypersurface type\/}, because a real hypersurface $M$
of a complex manifold $N$ carries such CR structure. \par
For a CR structure of hypersurface type,
the distribution $\Cal D$ can be locally described  as the kernel of
a  1-form $\theta$. Such  form $\theta$ determines
 an Hermitian metric
$${\Cal L}^\theta_q\: \D_q\times\D_q \to \R$$
 by the formula
$${\Cal L}^\theta( v, w) = (d\theta)(v, Jw) $$
for  any $v,w \in \Cal D$. This form is called the {\it Levi form
 of $(M,\D, J)$
associated with the form $\theta$\/}.
Notice  that the 1-form $\theta$ is defined up to  multiplication
by a function $f$ everywhere different from zero and
that ${\Cal L}^{f\theta} = f{\Cal L}^{\theta}$. In particular, the conformal
class
of a Levi form depends only on the CR structure.
\par
A  CR structure $(\D, J)$ of hypersurface type is called
{\it Levi non-degenerate} if it has non-degenerate Levi form or, in other
words, {\it if $\D$ is a contact distribution\/}. In this case we will call
$\D$ {\it the contact distribution underlying the non-degenerate contact
structure $(\D, J)$\/}.
\par
\medskip
Let $(M, \D)$ and $(M', \D')$ be two contact manifolds with contact
distributions $\D$ and $\D'$, respectively.
  A smooth map
$\varphi\: M \to M'$ is called {\it contact map\/} if
$\varphi_*(\D) \subset \D'$.\par
  A smooth map
$\varphi\: M \to M'$ of a  CR manifold $(M, \D, J)$
into some other  CR manifold  $(M', \D', J')$ is called
 {\it  holomorphic map  \/} or  {\it CR map} if \par
a) $\varphi_*(\D) \subset \D'$;\par
b) $\varphi_*(Jv) = J'\varphi_*(v)$ for all $v\in \D$.
\par
 In particular, we  define a {\it CR  transformation} of a CR
manifold $(M, \D, J)$ (resp. {\it contact transformation\/})
as a transformation $\varphi: M\to M$ such that
$\varphi$ and $\varphi^{-1}$ are both CR maps (resp. {\it contact map\/}).
It is known that
the group $Aut(M,\D, J)$ of all CR transformations of a Levi non-degenerate
CR manifold is a Lie group.
\par
\medskip
If the opposite is not stated, by CR manifold we will mean
{\it a simply connected Levi non-degenerate CR manifold\/}.
\medskip
We will also adopt the following notation. The symbols
$$ A =A(M, \D, J)  ,\qquad
A^{ss} = A^{ss}(M, \D, J) $$
 denote a   maximal connected compact subgroup and a maximal
connected compact semisimple subgroup
of $\text{Aut}(M, \D, J)$, respectively.  Recall that any two
maximal connected compact subgroups (resp. maximal connected compact
semisimple subgroups) are
conjugated  by an element of  $\text{Aut}(M, \D, J)$ .\par
\medskip
The Lie algebra of a Lie group is always denoted by the corresponding
gothic letter. For a subset $B$ of a Lie group $G$ or of a Lie algebra $\frak g$,
we denote by $C_G(B)$ and $C_{\g}(B)$ its centralizer in $G$ and
$\g$, respectively.\par
$Z(G)$ and $Z(\g)$ denote the center of a Lie group $G$ and of a Lie
algebra $\g$, respectively.\par
For any compact Lie group $G$ and the corresponding Lie algebra $\g$,
the expressions $G = G^{ss}\cdot Z(G)$ and $\g = \g^{ss} + Z(\g)$ denote the
decomposition into semisimple part plus center of $G$ and $\g$, respectively.\par
\medskip
For a compact Lie group $G$,  we  will  denote by $\B$ an
 $\operatorname{Ad}(G)$-invariant scalar
product   on the Lie algebra  $\g$. For example, if $G$ is simple,
$\B$ is a multiple of the Cartan-Killing form of $\g$. Throughout the paper, any orthogonal
decomposition of the Lie algebra $\g$ has to be understood
as orthogonal with respect to the inner product $\B$.\par

\smallskip
 By a homogeneous manifold $M= G/L$,  we mean a simply connected homogeneous manifold
of a connected  Lie group $G$, which acts almost effectively on $M$ (i.e. with
discrete kernel of non-effectivity).
It follows that  the stability subgroup $L$ is connected.

\bigskip
\subsubhead 2.2 First properties of homogeneous CR manifolds of compact Lie groups
\endsubsubhead
\medskip
 In this section, we recall some elementary facts about infinitesimal description of
contact and CR homogeneous manifolds.\par
\smallskip
Let $M = G/L$ be a homogeneous manifold of a connected compact Lie group
$G$  and $ \g = \l + \l^{\perp}$   the associated $\B$-orthogonal decomposition
of $\g$. Recall that $\l^\perp$ is naturally
identified with $T_{eL} G/L$. \par
Now, let  $\D \subset TM = T(G/L)$ be a $G$-invariant
contact distribution and $\theta$  a $G$-invariant contact form, i.e.
a 1-form such that $\theta(X) = 0$ for any $X\in \D$.
The {\it Reeb  field  associated with $\theta$\/} is the unique vector
field $\xi^\theta$ on $M = G/L$ such that
$$\theta_p(\xi^\theta) = 1\ ,\qquad d\theta_p(\xi^\theta, *)|_{\D} = 0\ ,\qquad \forall
p\in M\ .$$
By identifying $T_{eL} M$ with $\l^\perp$, we get a natural
$\operatorname{Ad}_L$-invariant decomposition
$$\l^\perp = \R Z + \m$$
where $\R Z$ is the 1-dimensional subspace corresponding to $\R \xi^\theta_{eL}
\subset T_{eL} M$ and $\m$ is the codimension one subspace corresponding
to the subspace $\D_{eL} \subset T_{eL} M$.\par
One can check that the decomposition
$$\g = \l + \R Z + \m\tag 2.1$$
is $\B$-orthogonal and that
the element
$Z \in \l^{\perp}$, defined up to a scaling, generates a  closed 1-parametric
 subgroup of $G$ and has a centralizer $\k = C_\g(Z)$, which is equal to
$\l \oplus \R Z$ (see \cite{AS}). \par
Any element $Z \in \l^{\perp}$, which generates a  closed 1-parametric
 subgroup and such that $C_\g(Z) = \l + \R Z$, is called {\it contact element
of $G/L$\/}. The formula (2.1) establishes a 1-1 correspondence between contact elements $Z$
up to scaling and $G$-invariant contact distributions $\D$,  $\D_{eL} = \m$.
\par
\medskip
The adjoint orbit
 $$F_Z = \operatorname{Ad}_G Z = G/K\ , \qquad K =C_G(Z)$$
is a flag manifold (i.e. a homogeneous manifold of a compact
semisimple Lie group $G$, which is $G$-isomorphic to an adjoint
orbit of $G$) and it is  called  {\it
the flag manifold $G$-associated with the
  contact manifold\/}
$(M=G/L,\D_Z)$. The $S^1$-fibration
$$  M=G/L \rightarrow F_Z = G/K $$
is called  {\it the structural $G$-fibration\/} of $(M = G/L, \D_Z)$.\par
\medskip
The reader should be aware that,  if $\B$ and $\B'$ are two
$\operatorname{Ad}_G$-invariant scalar products on $\g$, then the corresponding
contact  elements $Z\in \l^\perp $ and $Z'\in \l^\perp{}'$
associated with a given Reeb vector field $\xi^\theta$,
are in general different. Nevertheless, they verify
 $Z' = Z \mod Z(\l)$. Therefore
$$\k = C_\g(Z) = \l + \R Z = \l + \R Z' = C_\g(Z') = \k'\ .$$
This means that the $G$-associated flag manifold  $F_Z = F_Z'$
is independent
on the choice of the invariant inner product $\B$.\par
\medskip
In case $G = A^{ss}(M, \D, J)$ is a maximal compact semisimple
group of CR transformations of a homogeneous
CR manifold,
we call $F_Z$  {\it the flag manifold naturally associated with $(M, \D, J)$\/} and the
structural $G$-fibration $\pi: M \to F_Z$ is
 called {\it the structural fibration associated with $\D$\/}.\par
\smallskip
Now, let us
choose  a flag manifold  $F = G/K$. In \cite{AS} it was proved that
any homogeneous contact manifold $(M=G/L, \D_Z)$,  which have $F$
  as  $G$-associated
flag manifold, is obtained as follows. \par
An element $Z \in Z(\k)$ is
called $\k$-regular if:
\roster
\item"a)" $Z$ generate a closed 1-parametric subgroup of $G$;
\item"b)" $C_\g(Z) = \k$.
\endroster
If $Z$ is a $\k$-regular element,  the subalgebra
$\l_Z = \k \cap (Z)^\perp$ generates a closed
connected subgroup $L_Z \subset G$.  Moreover:
\medskip
\proclaim{Proposition 2.1}[AS] Let $F = G/K$ be a flag manifold
of a connected, compact semisimple Lie group $G$. There exists
a natural 1-1 correspondence
$$Z \quad \Leftrightarrow\quad (G/L_Z, \D_Z)$$
between $\k$-regular elements $Z \in \g$ (determined up to a scaling)
and  homogeneous contact manifolds $(G/L, \D)$,  which have
$F = G/K$ as
$G$-associated flag manifold.
\endproclaim
\bigskip
For a given  flag manifold $G/K$ and a $\k$-regular element
$Z \in \g$, we  say that $(G/L_Z, \D_Z)$ is
{\it the contact homogeneous manifold associated with
the pair $(G/K, Z)$\/}.\par
\par
\bigskip
Any $G$-invariant CR structure $(\D_Z, J)$ on $G/L$,
with underlying contact distribution $\D_Z$
and with associated   decomposition (2.1),
is uniquely associated with a
complex subspace $\m^{10}$ of $\m^\C$ such that
\roster
\item"i)" $\m^{10} \cap \overline{\m^{10}} = \{0\}$ and $ \m^\C=  \m^{10} + \overline{\m^{10}} $;
\item"ii)" $\l^\C + \m^{10}$ is a complex subalgebra of $\g^\C$.
\endroster
In fact, the decomposition $\m^\C = \m^{10} + \m^{01} = \m^{10} + \overline{\m^{10}}$
defines a complex structure $J: \m \to \m$ which corresponds to the complex structure
$J : \D \to \D$ of the distribution $\D$.\par
We call $\m^{10}$ {\it the holomorphic subspace $G$-associated with
$(\D,J)$\/}, since the associated complex sub-distribution
$T^{10} M \subset \D^\C$ is the eigenspace distribution of $J$ with
eigenvalue $i$.\par
\medskip
  Note that if we denote by  $\k = C_{\g}(Z)$ and if
the subspace  $\m^{10}$ is $\text{ad}_{\k}$-invariant,
 then $\m^{10}$ defines also an invariant complex structure $J_F$ on the flag
 manifold $F_Z =G/K$, since we may identify  $T_o^{10}F = \m^{10}$ where $o = eK$
(see \cite{AS} for more details). \par
\bigskip
\definition{Definition 2.2} Let $(M = G/L, \D, J)$ be a  homogeneous  CR manifold
of a connected
compact semisimple Lie group $G$ with contact element  $Z \in \g$
 and  holomorphic subspace  $\m^{10} \subset \g^\C$.\par
We say that the CR structure on $M$ is {\it G-standard\/} if
$$[Z, \m^{10}] \subset  \m^{10} \ .$$
If $G = A^{ss}(M, \D, J)$, a G-standard CR structure
will be called {\it standard\/}.
\enddefinition
\bigskip
 Note that  a CR structure $(\D,J)$ is $G$-standard if and only if the
holomorphic subspace $\m^{10}$ is
$\operatorname{Ad}_K$-invariant, where  $K = C_G(Z_\g)$.
In this case  we will denote  the invariant complex structure
defined by the subspace $\m^{10}$  on $F$ be  $J_F$, and we will consider
the associated
 flag manifold $F=G/K$  as a complex homogeneous
manifold with complex structure $J_F$.\par
Then the canonical fibration
$$ \pi: M =G/L \rightarrow F= G/K $$
is  a holomorphic fibration with respect to the CR structure $(\D,J)$ on $M$ and
the complex structure $J_F$ (which may be  considered  as a codimension zero
CR structure). This is a characteristic property of standard
CR structures (see  \cite{AS}).\par
\bigskip
\subsubhead 2.3 Compact homogeneous CR manifolds as homogeneous manifolds
 of compact semisimple  Lie groups
\endsubsubhead
\medskip
It is known that if $(M, \D, J)$ is   compact and homogeneous, then
any maximal compact subgroup $A(M, \D, I) \subset  \text{Aut}(M, \D, I)$ acts transitively
on $M$ (see \cite{AHR} and \cite{Sp}).
This result together with  Cor. 3.2 in \cite{AS} can be used
to prove the following more precise result.\par
\bigskip
\proclaim{Proposition 2.3} Let $(M, \D, J)$ be a compact homogeneous CR manifold.
Then  a maximal connected compact semisimple subgroup $A^{ss}(M, \D,J)$  acts
transitively on $M$.
\endproclaim
\demo{Proof} By the above remarks, we may represent the manifold $M$ as
the homogeneous manifold
 $M = G/L$, where $G = A(M, \D,J)$ is a maximal connected compact subgroup
of  $\text{Aut}(M, \D, J)$.
By Corollary  3.2  in  \cite{AS}, if $G$ is not semisimple, $G = G^{ss} \cdot Z(G)$
where the center $Z(G)$ has dimension one.
Moreover,  $\dim (\l \cap \g^{ss}) = \dim \l - 1$.
Hence, the semisimple part $G^{ss}$ acts transitively on $M$
because the $G^{ss}$-orbit of $o = eL$ has dimension
$\dim G^{ss}\cdot o = \dim \g^{ss} - \dim (\l\cap \g^{ss})
= \dim G/L$. \qed
\enddemo
\bigskip
\bigskip
\subhead 3. Characterization of G-standard CR manifolds
\endsubhead
\bigskip
Let $(M = G/L, \D, J)$ be  a  homogeneous  CR manifold
of a connected compact semisimple Lie group $G$. We will prove that the
 property of being $G$-standard does not depend on the
group $G$.\par
Let $A \subset \text{Aut}(M, \D)$ be a connected semisimple
group of contact automorphisms,
which contains $G$. Then $M = A/B = G/L$
 where $L = G\cap B$ and
 we have the orthogonal decomposition
$$\g = \l + \R Z_\g + \m_\g\ ,\qquad \a = \b + \R Z_\a + \m_\a\  $$
associated with the  invariant contact structure
$\D$, where $Z_\g$ and $Z_\a$ are the corresponding
contact elements. \par
If $A$ is  a group of CR transformations, then we  denote by
$\m^{10}_\g \subset \m^\C_\g$
and $\m^{10}_\a \subset \m_\a^\C$
the
 holomorphic
subspaces  of the CR manifolds $(G/L,\D, J)$ and $(A/B, \D, J)$.
\par
\bigskip
\proclaim{Lemma 3.1} Let $(M = G/L, \D)$ be a  homogeneous  contact manifold
of a compact connected semisimple Lie group $G$ and
$A\supset G$  a compact connected semisimple subgroup
of $\text{Aut}(M, \D)$ so that $M = A/B = G/L$ and $L = G\cap B$.  Then:
\roster
\item if $\pi: \a \to \m_\a$ is the projection parallel to $\b + \R Z_\a$,
then $\pi|_{\m_\g} : \m_\g \to \m_\a$ is an isomorphism;
\item $C_\g(Z_\a) = C_\g(Z_\g) = \l + \R Z_\g$.
\endroster
Moreover, if  $(\D, J)$ is a $G$-invariant CR structure on $G/L$
 with underlying  contact structure $\D$ and if
$A \subset \text{Aut}(M , \D, J)$, then:
\roster
\item"(3)" the isomorphism  $\pi|_{\m_\g} : \m_\g \to \m_\a$
 commutes
with the complex structures on $\m_\g$ and $\m_\a$ induced by the CR structure
$(\D, J)$ and
$$\l^\C + \m^{10}_\g \subset \b^\C + \m^{10}_\a\ ;$$
\item"(4)" $[Z_\g, \m^{10}_\g] \subset \l^\C + \m^{10}_\g $ if and only if
$[Z_\a, \m^{10}_\a] \subset \b^\C + \m^{10}_\a$.
\endroster
\endproclaim
\demo{Proof} We denote by $\B$ an invariant scalar product on $\a$ and
we may assume that the invariant scalar product on $\g$ is the restriction
of $\B$ to $\g\subset \a$. \par
Then (1)  follows immediately
from the fact that $(\l + \m_\g)/\b \simeq
 (\b + \m_\a) / \b \simeq \D_{o}$ (here $o = e L = eB$).\par
\smallskip
 To show (2), note that  $C_\a(Z_\a) = \b + \R Z_\a$ and
$C_\g(Z_\g) = \l + \R Z_\g$. Hence, it is sufficient to check that
$$(\b + \R Z_\a) \cap \g = \l + \R Z_\g\ .$$
But this follows immediately from the fact that
$X\in \b + \R Z_\a$ (resp. $X\in \l + \R Z_\g$) if and only if
the value at $o$ of the corresponding vector field  $\hat X$ is proportional
with the Reeb vector $\xi^\theta_{o} \in T_oM$ at $o$.\par
\smallskip
For (3), observe that  the first claim follows directly from the definitions. The second claim
is proved by the  fact that $\m^{10}_\g =
\m^{10}_\a \mod \b^\C$.\par
\smallskip
Let us now prove (4). From the  proof of (2), after  rescaling $Z_\g$,  we may assume that
$$Z_\a = Z_\g + W$$
for some $W \in \b$.  Moreover, by (1) and (3),
 any element $X^{10}_\a \in \m^{10}_\a$ can be written as
$$X^{10}_\a = X^{10}_\g + X_\b$$
where $X^{10}_\g  \in \m^{10}_\g$ and $X_\b \in \b^\C$.
\par
Assume that  $[Z_\g, \m^{10}_\g] \subset \l^\C + \m^{10}_\g$. Then for any
$X^{10}_\a =  X^{10}_\g + X_\b \in \m^{10}_\a$ we  have
$$[Z_\a, X^{10}_\a ] = [Z_\g + W, X^{10}_\a] \cong [Z_\g, X^{10}_\a]
(\mod (\b^\C + \m_\a^{10})) \cong $$
$$ \cong [Z_\g, X^{10}_\g] + [Z_\g, X_\b] \ (\mod (\b^\C + \m_\a^{10})) =$$
$$  = [Z_\g, X^{10}_\g] + [Z_\a, X_\b]  - [W, X_\b] \ \mod (\b^\C + \m_\a^{10})
=$$
$$ = [Z_\g, X^{10}_\g] \mod (\b^\C + \m_\a^{10})  = $$
$$ =
0\  \mod (\b^\C + \m_\a^{10}) $$
by assumption and  (3).\par
Conversely, assume that  $[Z_\a, \m^{10}_\a] \subset \b^\C + \m^{10}_\a$. Then
$$[Z_\g, \m^{10}_\g] \subset [Z_\g, \b^\C + \m^{10}_\a]\cap \g^\C
= [Z_\a - W, \b^\C + \m^{10}_\a] \cap \g^\C \subset
$$
$$\subset (\b^\C + \m^{10}_\a) \cap \g^\C = \l^\C + \m^{10}_\g$$
by assumption and  (3).
\qed
\enddemo
\bigskip
Lemma 3.1 implies  the following proposition.\par
\medskip
\proclaim{Proposition 3.2} Let $(M = G/L, \D, J)$ be a homogeneous contact manifold
and  $A \subset  \text{Aut}(M, \D)$  a connected compact semisimple group of
contact transformations  which contains $G$,  so that   $M = A/B = G/L$.
Denote by $Z_\g$ and $Z_\a$  contact elements in $\g$ and $\a$ associated
with the contact distribution $\D$. Then:
\roster
\item"a)" the $G$-associated flag manifold $F_{Z_\g} = G/K = G/C_G(Z_\g)$ is equivalent (as $G$-manifold)
to
the $A$-associated flag manifold $F_{Z_\a} = A/C = A/C_A(Z_\a)$.\par
\item"b)" if  $A$ is  a group of CR transformations, then
$(M = G/L= A/B, \D, J)$ is a G-standard CR manifold if
and only
if it is an A-standard CR manifold.
\endroster
\endproclaim
\demo{Proof} a) It follows from the fact that the orbit
$G\cdot o \subset A/C_A(Z_\a)$ coincides with $A/C_A(Z_\a) = F_{Z_\a}$
and it is equal to
$G/C_G(Z_\a) = G/C_G(Z_\g)$, by Lemma 3.1 (2).\par
b)  follows directly from definitions and Lemma 3.1 (4).\qed
\enddemo
\bigskip
We have the following  corollary, which
is the main result of this section.\par
\medskip
\proclaim{Corollary 3.3} For any connected semisimple Lie group
$G$ acting transitively on a compact CR manifold $(M, \D, J)$, the
$G$-associated flag manifold $G/K$  is $G$-equivariantly
diffeomorphic to  the naturally associated flag manifold \par
\noindent
 $F = A^{ss}(M,\Cal D , J )/B$ naturally associated with
$(M, \Cal D, J)$.
\par Moreover $(M, \D, J)$ is
 G-standard if and only if it is standard.
\endproclaim
\bigskip
\bigskip
\subhead 4.  Maximal compact semisimple groups of automorphisms
\endsubhead
\bigskip
\subsubhead 4.1 Inclusion relations between transitive groups of
 transformations of a flag manifold
\endsubsubhead\par
\medskip
First of all we quote the following result by A. L. Onishchik, which describes the
inclusion relations between compact semisimple transitive groups of holomorphic
transformations of a flag manifold.\par
\bigskip
\proclaim{Theorem 4.1} [On] Let $F = G/K$ be a flag manifold of a compact
semisimple Lie group $G = G_1\times \dots \times G_p$, where
 $G_i$ are the simple factors of $G$, and  let
$$F = G/K  = G_1/K_1\times \dots \times G_p/K_p$$
 be the corresponding decomposition of the flag manifold $F$.\par
Let $A$ be a  compact semisimple Lie group
of transformations
of $F$, which contains $G$ and preserves some complex structure $J$ on $F$.
Then $A$ is of the form
$$A = A_1\times \dots \times A_p$$
where each $A_i \supseteq G_i$ acts  only on $F_i = G_i/K_i = A_i/C_i$ and
either $A_i = G_i$ or the pair  $(A_i/C_i, G_i/K_i)$
is one of those contained in the following table: \par
\medskip
\centerline{\bf Table 1}
\medskip
\moveright2.5cm
\vbox{\offinterlineskip
\halign {\strut\vrule\hfil\ $#$\ \hfil
 &\vrule\hfil\ $#$\
\hfil&\vrule\hfil\  $#$\
\hfil&\vrule\hfil\ $#$\
\hfil\vrule\cr
\noalign{\hrule}
\phantom{\frac{\frac{1}{1}}{\frac{1}{1}}}
n^o
&
F = A/C = G/K
&
 A/C
&
G/K
\cr \noalign{\hrule}
I
&
{}^{\phantom{e^{\frac{A^A}{A_A}}}}\C P^{2 \ell-1}\
\ell >1{}_{\phantom{e^{\frac{A^A}{A_A}}}}
&
\frac{SU_{2\ell}}{U_{2 \ell - 1}}
&
\frac{Sp_{\ell}}{Sp_{\ell - 1} \cdot T^1}
\cr \noalign{\hrule}
II
&
{}^{\phantom{e^{\frac{A^A}{A_A}}}}
Gr_2(\R^7)
{}_{\phantom{e^{\frac{A^A}{A_A}}}}
&
\frac{SO_{7}}{SO_{5}\cdot SO_2}
&
\frac{G_2}{U_{2}}
\cr \noalign{\hrule}
III
&
{}^{\phantom{e^{\frac{A^A}{A_A}}}}
Com(\R^{2\ell + 2})\ \ell >2
{}_{\phantom{e^{\frac{A^A}{A_A}}}}
&
\frac{SO_{2\ell + 2}}{U_{\ell+1}}
&
\frac{SO_{2\ell + 1}}{U_{\ell}}
\cr \noalign{\hrule}
}}
\medskip
In  row II, the manifold $G_2/U_2$ is equal to $\operatorname{Ad}_{G_2}(H_\alpha) =
G_2/C_{G_2}(H_\alpha)$,
where $H_\alpha$ is the dual vector of a long root of $G_2$.
\endproclaim
\medskip
In the following, we will call {\it Onishchik pair\/} any pair  $(A/C, G/K)$
of homogeneous spaces from a row of Table 1.\par
\bigskip
As a direct corollary of Onishchik's result,
we get the following theorem.\par
\medskip
\proclaim{Theorem 4.2}
 Let $F = G/K$ be a flag manifold  and
$A$  a connected compact group of
transformations of $F$ which contains $G$ and
preserves some complex structure on $F$.
Then any $G$-invariant complex structure on $F$
is $A$-invariant.
\endproclaim
\demo{Proof} By Theorem 4.1, it is sufficient
to consider only the three cases of Table 1.
In all such  cases, $A/C$ is an irreducible  Hermitian symmetric
manifold and hence it admits a unique (up to a sign) invariant complex
structure. It is also easy to check that also  $M = G/K$, where
$M$ is  $Sp_\ell/Sp_{\ell-1}\cdot T^1$,
$SO_{2\ell + 1}/U_\ell$ or
  $G_2/U_2 = G_2/C_{G_2}(H_\alpha)$, with $\alpha$ long root, admits
only one  (up to a sign)
invariant complex structure. In fact, any
invariant complex structures  on a flag
manifold $G/K$ (considered up to conjugation) corresponds to a black-white
Dynkin diagram of the Lie group $G$, where the subdiagram formed by the white nodes
is equal to the Dynkin diagram of the semisimple part $K'$ of $K$
(see e.g. \cite{Al}, \cite{AP}); for the three manifolds  $M = G/K$ above,
there is only  one black-white
Dynkin diagram.
\qed
\enddemo
\remark{Remark 4.3} From Theorem 4.2, it follows immediately that if
$F = G/K = G_1/K_2 \times \dots G_p/K_p$ is the decomposition into irreducible
factors of a flag manifold with invariant complex structure $(F, J_F)$,
then the maximal group $A$ of holomorphic
transformations  is $A = A_1 \times \dots \times A_p$,
where $A_i = \tilde A$ if there exists an Onishchik pair of the form $(\tilde A/\tilde C, G_i/K_i)$
and it is $A_i = G_i$ otherwise.
\endremark
\bigskip
\subsubhead 4.2 Homogeneous manifolds with non-standard CR structures
\endsubsubhead
\par
\medskip
By Corollary 3.3, a homogeneous CR manifold
$(M = G/L, \D, J)$ of a compact semisimple Lie group $G$ is not $G$-standard
 if and only if it is non-standard. We recall the classification of such
non-standard CR manifolds
 in the following Theorem (see \cite{AS}). \par
\bigskip
\proclaim{Theorem 4.4} [AS]
Let  $(M = G/L, \D_Z, J)$ be a simply connected
non-standard CR manifold of a compact connected Lie group $G$
 and  $F = G/K$ the associated  flag manifold. Then
the triple $(G, L, K)$ is one of those given in  Table 2.\par
\medskip
\moveright 1 cm
\vbox{\offinterlineskip
\halign {\strut\vrule\hfil\ $#$\ \hfil
 &\vrule\hfil\ $#$\
\hfil&\vrule\hfil\ $#$\
\hfil&\vrule\hfil\  $#$
\hfil
\vrule\cr
\noalign{\hrule} n^o &
\phantom{\frac{\frac{1}{1}}{\frac{1}{1}}}G
\ \ &
L
&
K
\cr \noalign{\hrule}
1 &
SU_2\times SU'_2
&
\phantom{\frac{\frac{1}{1}}{\frac{1}{1}}}
T^1
&
T^1\times T^1{}'
\cr \noalign{\hrule}
2
&
Spin_7
&
\phantom{\frac{\frac{1}{1}}{\frac{1}{1}}}
SU_3
&
T^1\cdot SU_3
\cr \noalign{\hrule}
3
&
F_4
&
\phantom{\frac{\frac{1}{1}}{\frac{1}{1}}}
Spin_7
&
T^1\cdot SO_7
\cr \noalign{\hrule}
4 &
SU_{2}
&
\phantom{\frac{\frac{1}{1}}{\frac{1}{1}}}
\{ e\}
&
T^1
\cr \noalign{\hrule}
5 &
SO_{2n+1} \smallmatrix n >1 \endsmallmatrix
&
\phantom{\frac{\frac{1}{1}}{\frac{1}{1}}}
SO_{2n-1}
&
T^1\cdot SO_{2n-1}
\cr \noalign{\hrule}
6
&
SO_{2n} \smallmatrix n >2 \endsmallmatrix
&
\phantom{\frac{\frac{1}{1}}{\frac{1}{1}}}
SO_{2n-2}
&
T^1\cdot SO_{2n-2}
\cr \noalign{\hrule}
7
&
Sp_n
&
\phantom{\frac{\frac{1}{1}}{\frac{1}{1}}}
Sp_1\cdot Sp_{n-2}
&
T^1\cdot Sp_1\cdot Sp_{n-2}
\cr \noalign{\hrule}
8
&
SU_{n}\
&
\phantom{\frac{\frac{1}{1}}{\frac{1}{1}}}
T^1 \cdot SU_{n-2}
&
T^1 \cdot U_{n-2}
\cr \noalign{\hrule}
9
&
\underset{\phantom{B}}\to{
\overset{\phantom{A}}\to{\matrix SU_p\times SU'_q\\
\smallmatrix p + q > 4\endsmallmatrix
\endmatrix}}
&
\phantom{\frac{\frac{1}{1}}{\frac{1}{1}}}
T^1\cdot U_{p-2} \cdot U'_{q-2}
&
(T^1\cdot U_{p-2}) \cdot (T^1{}' \cdot U'_{q-2})
\cr \noalign{\hrule}
10
&
SU_n\ \smallmatrix n >4 \endsmallmatrix
&
\phantom{\frac{\frac{1}{1}}{\frac{1}{1}}}
T^1\cdot (SU_2 \times SU_2) \cdot SU_{n-4}
&
T^1\cdot (SU_2 \times SU_2) \cdot U_{n-4}
\cr \noalign{\hrule}
11
&
SO_{10}
&
\phantom{\frac{\frac{1}{1}}{\frac{1}{1}}}
T^1\cdot SO_6
&
T^2\cdot SO_6
\cr \noalign{\hrule}
12
&
E_6
&
\phantom{\frac{\frac{1}{1}}{\frac{1}{1}}}
T^1\cdot SO_8
&
T^2\cdot SO_8
\cr \noalign{\hrule}
}}
\endproclaim
\noindent
\centerline{\bf Table 2}
\bigskip
\bigskip
\subsubhead 4.3 Maximal compact semisimple groups of automorphisms
of  non-standard CR manifolds
\endsubsubhead
\medskip
Using Onishchik's theorem (Theorem 4.1), Proposition 3.2 and Theorem 4.4,
 we can now describe the maximal compact semisimple group of CR transformations
for any compact  homogeneous non-standard CR manifold. The answer
 is quite simple.
\par
\bigskip
\proclaim{Theorem 4.5} Let $(M=G/L,\D, J)$ be a non-standard homogeneous
CR manifold of a semisimple connected compact Lie group $G$. \par
\roster
\item If $M  = G/L\neq  \frac{\text{Spin}_7}{SU_3}$,
then $G$ is
a maximal connected semisimple compact Lie group of automorphisms.
\item If $M =  \frac{\text{Spin}_7}{SU_3}$,
then $A = SO_8$ is a maximal connected semisimple compact Lie group of automorphisms
of $M$ and  $M = \frac{SO_8}{SO_6} = \frac{\text{Spin}_7}
{SU_3}$.
\endroster
\endproclaim
\demo{Proof} Assume that there exists a compact semisimple group $A$ of
automorphism of $M$ which properly contains $G$, so that $M = G/L = A/B$,
with $L = G \cap B$. By Proposition 3.2 (a), the $A$-associated
flag manifold $F = A/C$ is equivalent to the $G$-associated flag manifold $G/K$
and hence, by  Theorem 4.1, one of the factors of $F = A/C$
must be a member of an Onishchik pair.\par
At the same time, $(M = A/B, \D, J)$ is non-standard and hence
 the groups $A$, $B$ and  $C$ must be in a row
 of Table 2.\par
Comparing Table 1 and Table 2, we find
 that there  are only two  possibilities for $F = A/C$, that is
 $F = \frac{SO_{7}}{SO_{5}\cdot SO_2}$ or
$F = \frac{SO_{8}}{U_4} = \frac{SO_8}{T^1 \cdot SO_6}$. It follows that
$G/K = \frac{G_2}{U_{2}}$ or $G/K = \frac{SO_7}{U_3} = \frac{Spin_7}{U_3}$.
Since  by Table 2 there is no non-standard homogeneous CR manifold
 with $G = G_2$, it follows that the first case is  impossible and hence that (1)
is true.
\par
In order to prove (2),
notice that for any non-standard homogeneous
CR manifold $(M = \frac{\text{Spin}_7}{SU_3}, \D, J)$,
the anticanonical map (see definition in \cite{AHR}; see also \cite{AS}, \S 4.4)
determines a CR equivalence between $M$ and a real hypersurface
$$\frac{\text{Spin}_7}{SU_3} =
\text{Spin}_7 \cdot v \subset
 T S^7\ ,\qquad \text{for some} \ 0 \neq v\in TS^7 = T \left(\frac{\text{Spin}_7}{G_2}\right)\ ,$$
where we identify $T S^7$ with the  complex homogeneous space
$T S^7 = T \left(\frac{SO_8}{SO_7}\right) = \frac{SO_8(\C)}{SO_7(\C)}$. On the other hand,
it is clear that for any $0 \neq v\in TS^7$ we have
$$\frac{\text{Spin}_7}{SU_3} = \text{Spin}_7 \cdot v = SO_8 \cdot v = \frac{SO_8}{SO_6}\ .$$
Since  $SO_8$ acts on $T S^7 = \frac{SO_8(\C)}{SO_7(\C)}$ as
a group of biholomorphisms, it follows that it
acts on $M = \frac{\text{Spin}_7}{SU_3} \simeq \text{Spin}_7 \cdot v$ as a transitive group of
CR transformations which properly contains $\text{Spin}_7$.
By the proof of (1),  it follows also that  $SO_8$  is
a maximal compact semisimple Lie group of CR transformations of $M$.
\qed
\enddemo
\bigskip
\subsubhead 4.4 Maximal compact semisimple groups of automorphisms
of  standard CR manifolds
\endsubsubhead
\medskip
It remains to describe the maximal connected,
compact semisimple group
of CR transformations  of a standard CR manifold.  \par
For this, we first need the following
Lemma, where we denote by $\B_\a$ (resp. $\B_\g$)
an invariant scalar product on $\a$ (resp. $\g$).\par
\bigskip
\proclaim{Lemma 4.6} Let $F = A/C = G/K$  be a flag manifold
with $A \supsetneq G$, and  $(M = G/L, \D)$ and
$(M' = A/B, \D')$ be two  homogeneous contact manifolds, having
$F$ as associated flag manifold. Denote by
$Z_\g\in \g$ and $Z_\a\in \a$  the contact elements
associated with $\D$ and $\D'$, respectively.\par
Then the following two properties
are equivalent:
\roster
\item"i)"  $G$  acts transitively on  $M' = A/B$ and the homogeneous manifold
 $M' = A/B = G/L'$   is $G$-equivalent to
$M = G/L$;
\item"ii)"
 $$Z(\k) \cap Z_\g^{\perp_\g} = Z(\k) \cap Z_\a^{\perp_\a}\ ,\tag4.1$$
 where
$Z_\g^{\perp_\g} = \{ \ X \in \g\ : \B_\g( X, Z_\g) = 0\ \}$ and
$Z_\a^{\perp_\a} = \{ \ X \in \a\ : \B_\a( X, Z_\g) = 0\ \}$.
\endroster
Furthermore, if this is the case, then under the
identification $M = A/B = G/L$, where $L = B\cap G$,
the invariant contact structures $\D$ and
$\D'$ coincide.
\endproclaim
\demo{Proof} Let $F = G/K = G_1/K_1 \times \dots G_p/K_p$ be the decomposition of
$F$ into irreducible factors. By Theorem 4.1,
$$F = A/C = A_1/C_1 \times \dots \times A_p/G_p$$
where for each $i$, either $A_i/C_i = G_i/K_i$ or  $(A_i/C_i, G_i/K_i)$
 is an   Onishchik pair.\par
\medskip
Note that
$G$  acts transitively on  $M' = A/B = G/L'$ and that $G/L'$ is equivalent
to $G/L$ if and only if $\l = \l' = \b \cap \g$.
On the other hand, since
 $\l = \k \cap Z_\g^{\perp_\g} = \k^{ss} + Z(\k) \cap Z_\g^{\perp_\g}$
and $\b = \c \cap Z_\a^{\perp_\a} = \c^{ss} + Z(\c) \cap Z_\a^{\perp_\a}$,
we also have that
$$\b \cap \g = (\c \cap \g)\cap Z_\a^{\perp_\a} = \k \cap Z_\a^{\perp_\a} =
\k^{ss} + Z(\k) \cap Z_\a^{\perp_\a}\ .$$
Hence  $\l = \l' = \b \cap \g$ if and only if (4.1) holds. \par
Assume now that (4.1) holds and that $\D$ and
$\D'$ are not equal. Then  $M = G/L = A/B$ admits two
invariant contact structure and hence it
is a special contact manifold according to Def. 3.4
in \cite{AS}. Then,  by Thm. 3.6 of \cite{AS},
$G$ and $A$ are both simple.  Theorem 4.1 implies that
$(A/C, G/K)$ is an
Onishchik pair and, in particular, that the center $Z(\k)$ is
one dimensional. By Proposition 2.1, this implies that
there exists only  one homogeneous contact manifold $(G/L, \D)$
with associated flag manifold $F = G/K$; this is in contradiction with
the fact that $(G/L, \D)$ and $(G/L, \D')$ are two distinct invariant contact
homogeneous manifolds with associated $F$. \qed
\enddemo
\bigskip
Now we can prove the main theorem of this subsection.\par
\bigskip
\proclaim{Theorem 4.7} Let $(M = G/L, \D, J)$ be a standard homogeneous
CR manifold of a semisimple connected compact Lie group $G$. Let
also  $F = G/C_G(Z) = G/K$ be
the associated  flag manifold and $J_F$ the  complex structure
induced by the projection $\pi: G/L \to G/K$. \par
Then the maximal connected compact semisimple group $A^{ss}(M, \D,J)$
which contains $G$ acts on  $(F = G/K, J_F)$ as a
maximal compact group  of holomorphic transformations.
\endproclaim
\demo{Proof}  By Proposition  3.2, any
connected compact semisimple group,
which contains $G$, acts naturally on $F$ as a group of holomorphic
transformations.  Hence it is sufficient to prove that
a maximal connected semisimple group $A$ of holomorphic transformations of $(F, J_F)$,
which contains $G$, acts on $M = G/L$ as a group of CR transformations. \par
We first show that $A$ acts on $(M = G/L, \D)$ as a group of
contact transformations.\par
\medskip
We may assume that $\B = \B_\g$ is  the Cartan-Killing form of $\g$.\par
 Let
$$F = G/K = G_1/K_1 \times \dots \times G_p/K_p\ ,\quad Z_\g = Z_1 + \dots + Z_p\ ,\quad Z_i \in \g_i$$
 be the decomposition of $F$ into irreducible
factors and
  the associated decomposition of the contact element.
We may also
assume  that $G$ is simply connected and that
 $\exp(Z_\g) = \exp(Z_1)\cdot \dots \cdot \exp(Z_p) = e$  where
$e$ is the identity element of $G$. Clearly, this implies that
$\exp(Z_i) = e\in G_i$ for any $i$.\par
Recall that, since $Z_\g \in Z(\k)$ and
$\l = \k \cap (\R Z)^\perp$, $\g$ has the following
$\B_\g$-orthogonal decomposition:
$$\g = \l + \R Z  + \m = [(\k^{ss}_1 + \k^{ss}_2 + \dots +  \k^{ss}_p) + Z(\k) \cap (\R Z_\g)^\perp] +
\R Z_\g + (\m_1 + \dots + \m_p)$$
with   $\m_i = \m \cap \g_i$. \par
\medskip
By Remark 4.3,  $F$  can be decomposed into
$$F = G/K = A/C = A_1/C_1 \times \dots \times A_p/C_p$$
where either $A_i/C_i = G_i/K_i$ or $(A_i/C_i, G_i/K_i)$ is an Onishchik pair.
In case $(A_i/C_i, G_i/K_i)$ is an Onishchik pair,
 $\dim Z(\c_i) = \dim Z(\k_i) = 1$.  We
choose  a generator $E^{\c_i}$ for $Z(\c_i)$ and
a generator   $E^{\k_i}$ for $Z(\k_i)$, which verify
the following property: for any $X\in \R E^{\c_i}$ and any
$Y \in \R E^{\k_i}$, $\exp(X) = e$ and $\exp(Y) = e$ if and only if
$X\in \Z E^{\c_i}$ and $Y \in \Z E^{\k_i}$, respectively.
\medskip
We fix now an $\operatorname{Ad}_A$-invariant scalar product $\B_\a$ on $\a$.
For each simple algebra $\a_i$, we assume that $\B_\a|_{\a_i}$
is a multiple of the  Cartan-Killing form of $\a_i$
determined with the following rules: if
$A_i = G_i$,  we assume that $\B_\a|_{\a_i}$
is the Cartan-Killing form without rescaling (note that in this case, $\B_\a|_{\a_i} = \B_\g|_{\g_i}$);
 if $(A_i/C_i, G_i/K_i)$ is an Onishchik pair, we assume that $\B_\a|_{\a_i}$ is the
multiple of the Cartan-Killing form  which verifies
 $\B_\a(E^{\c_i}, E^{\c_i}) = -1$. \par
\medskip
We now consider  the  element $Z_\a = Z'_1 + \dots Z'_p \in Z(\c) \subset \a$
defined as follows: if $(A_i/C_i, G_i/K_i)$ is an Onishchik pair and $Z_i$ is of the form
$Z_i = \lambda_i E^{\k_i}\ ,$
  we set
$$Z'_i = \lambda_i\frac{\B_\g(E^{\k_i}, E^{\k_i})}{\B_\a(E^{\k_i}, E^{\c_i})} E^{\c_i}
= \frac{\B_\g(Z_i, Z_i)}{\B_\a(Z_i, E^{\c_i})} E^{\c_i}\ ;\tag 4.2$$
in all other cases,   we set $Z'_i = Z_i$. \par
 It is an immediate consequence of
definitions that  $\B_\a(Z_i, Z'_i) = \B_\g(Z_i, Z_i)$ for any $i$.\par
We claim that:
\roster
\item"a)" $Z_\a$ is a $\c$-regular element, that is $C_\a(Z_\a) = \c$ and $Z_\a$
generates a closed subgroup of $A$;
\item"b)" $Z(\k) \cap Z_\g^{\perp_\g} = Z(\k) \cap Z_\a^{\perp_\a}$
 where
$Z_\g^{\perp_\g} = \{ \ X \in \g\ : \B_\g( X, Z_\g) = 0\ \}$ and
$Z_\a^{\perp_\a} = \{ \ X \in \a\ : \B_\a( X, Z_\a) = 0\ \}$.
\endroster
Assume for the moment that a) and b) are true.  Then
by Proposition 2.1, there exists a unique   homogeneous contact manifold
$(A/B, \D_{Z_\a})$ with
contact structure  $\D_{Z_\a}$
 with  contact element $Z_\a$ and with $F = A/C = G/K$
as associated flag manifold; moreover, by Lemma 4.6,
$(M = G/L, \D)$ is equal to  $(A/B,\D_{Z_\a})$ and hence
$A$ acts on $M$ as a group of contact transformations, as we needed to prove.\par
To show  a), observe that
$$C_\a(Z_\a) = C_{\a_1}(Z'_1) + C_{\a_2}(Z'_2)  + \dots + C_{\a_p}(Z'_p)\ .$$
It is clear that   $C_{\a_i} = \c_i$ for any $i$ and hence that $C_\a(Z_\a) = \c$.
So, we only need  to check that $Z_\a$ generates a closed subgroup. For this,
 consider the following facts:
\roster
\item if $A_i = G_i$, then $\exp(Z'_i) = \exp(Z_i) = e \in G_i$;
\item if $(A_i/C_i, G_i/K_i)$ is an Onishchik pair,  then the elements
$E^{\c_i}$ and $E^{\k_i}$ are defined in such a way that $$\exp(Z_i) =
\exp(\lambda_i E^{\k_i}) = e\in G_i\ ,
\exp(\lambda'_i E^{\c_i}) = e\in A_i$$
imply $\lambda_i \in \Bbb
Z$ and $\lambda'_i \in \Bbb
Z$, respectively;
\item since $\exp(\R Z_i)$ is closed in $G_i$ for any $i$ and since
$\B_{\g}|_{\g_i}$ is the Cartan-Killing form of $\g_i$, we have that
$\B_\g(Z_i, Z_i) \in \Bbb Q$ for any $i$;
\item for any Onishchik pair $(A_i/C_i, G_i/K_i)$,
 consider the projection $\pi: C_i\to Z(C_i) = \exp(\R E^{\c_i}) = C/C^{ss}$;
the differential $\pi_*: \c_i \to Z(\c_i)$ is equal to the orthogonal projection
$$\pi_*: \c_i \longrightarrow Z(\c_i)\ ,\qquad \pi_*(X) = \B_\a(X, E^{\c_i}) E^{\c_i}\ ;$$
so
$$e = \pi(\exp(E^{\k_i})) = \exp(\pi_*(E^{\k_i})) = \exp(\B_\a(E^{\k_i}, E^{\c_i}) E^{\c_i})$$
and hence $\B_\a(E^{\k_i}, E^{\c_i})  \in \Bbb Z$;
from (2) this implies that
 $\B_\a(Z_i, E^{\c_i}) = \lambda_i \B_\a(E^{\k_i}, E^{\c_i})  \in \Bbb Z$.
\endroster
From (1),   (3) and (4) and  formula (4.2),
it follows immediately that there exists an integer $N$ such that
$\exp(N Z_\a) = e \in A$ and hence that  $\exp(\R Z_\a)$ is closed.
\smallskip
To prove b), consider  an
orthonormal  basis
$B_i = (Z_i, Y^i_2, \dots , Y^i_{p_i})$ for each abelian algebra $Z(\k_i)$, with first element
equal to $Z_i$. It follows that
an element
$$X = x^i Z_i + \sum_{j=2}^{p_i} y^j_i Y^i_j \in Z(\k)$$
is an element of  $Z(\k) \cap Z_\g^{\perp_\g}$ if and only if
$$x^1 \B_\g(Z_1, Z_1) + x^2 \B_\g(Z_2, Z_2) + \dots + x^p \B_\g(Z_p, Z_p) = 0\ .\tag 4.3$$
The same element is in  $Z(\k) \cap Z_\a^{\perp_\a}$ if and only if
$$x^1 \B_\a(Z_1, Z'_1) + x^2 \B_\a(Z_2, Z'_2) + \dots + x^p \B_\a(Z'_p, Z'_p) = 0\tag 4.4$$
Since $\B_\a(Z_i, Z'_i) = \B_\g(Z_i, Z_i)$ for any $i$, equations (4.3) and (4.4) coincide.\par
\medskip
It remains to check that  $A$ acts on $M = G/L$  as a group of CR transformations
of $(\D, J)$. Since the CR structure is standard and $A$ preserves the contact structure $\D$,
the above claim is  immediately proved by recalling
that $A$ is a group of holomorphic transformations for $(F, J_F)$.\qed
\enddemo
\bigskip
\bigskip
\subhead 5.  Maximal compact groups of automorphisms of a homogeneous
CR manifold
\endsubhead
\bigskip
In the previous section, we showed how to reconstruct a maximal compact
connected  semisimple
group of CR transformations $A^{ss}= A^{ss}(M, \D, J)$
of a  compact homogeneous CR manifold $(M = G/L, \D, J)$.
Now we want to show how to determine a maximal compact connected
group $A= A(M,\D,J)$  which contains $A^{ss}$.\par
\bigskip
\proclaim{Theorem 5.1} Let $(M, \D, J)$ be a simply connected  compact
 homogeneous CR manifold .
\roster
\item"(a)" If $(M, \D, J)$ is non-standard, then $A = A^{ss}$.
\item"(b)" If $(M, \D, J)$ is standard, then $A= A^{ss} \times T^1$.
\endroster
\endproclaim
\demo{Proof} (a) was proved in  [AS], Prop. 4.6.\par
(b) Assume now that $(M,\D,J)$ is standard and  identify $M$ with
 the  homogeneous manifold
$A^{ss}/B$, with $A^{ss} = A^{ss}(M,\D,J)$. Let
$$\a^{ss} = \b + \R Z + \m $$
be the associated decomposition. Without loss of generality
we may assume that the contact element  $Z$ is so that
$\exp(Z) = e$ and has $\B$-norm equal to $1$. Let also
 $\m^{10}\subset \m^{\Bbb C}$ be the
associated holomorphic subspace
(recall that it is $\text{ad}(\b + \Bbb R Z)$-invariant).
By Cor. 3.2 in \cite{AS},  the center $Z(A)$
of  $A = A(M,\D,J)$ has dimension 0 or 1. Hence
to prove (b) it is sufficient to define a CR action on $M$ of the group
$A = A^{ss}\times T^1$. This can be done by  constructing a homogeneous
CR manifold $(\tilde M = A/\tilde B ,\tilde {\D}, \tilde J)$ and then show
that it is
$A^{ss}$-diffeomorphic and CR-equivalent to  $(M=A^{ss}/B,\D,J)$.\par
Let $\tilde B$ be the connected subgroup of $A = A^{ss}\times T^1$
generated by the subalgebra
$$ \tilde {\b} =\b + \Bbb R (Z -\xi)  \subset \a = \a^{ss} \oplus \Bbb R \xi $$
where $\xi$ is a generator of $T^1$ such that $\exp(\xi) = e$. We may also
assume that the $\B$-norm of $\xi$ is equal to $1$.\par
One can   check easily that $\tilde Z = Z + \xi \in \tilde {\b}^{\perp} $
and that it is a contact element. Let
 $\tilde \D = \D_{\tilde Z}$ be the corresponding
contact structure on $\tilde M = A/\tilde B$. Note that
 the subspace $\m^{10}$
defines an invariant CR structure $(\tilde \D , \tilde J)$ on $\tilde M$ and
that the subgroup $A^{ss}$ of $A = A^{ss} \times T^1 $ acts transitively on
$\tilde M = A/\tilde B$ with stabilizer $B = A^{ss} \cap \tilde B $.
Therefore, $\tilde M = A/ \tilde B = A^{ss}/B = M$. It is quite simple to check that
the homogeneous CR manifold
$(\tilde M = A/\tilde B, \tilde \D, \tilde J)$
is CR equivalent
to  $(M = A^{ss}/B = \tilde M, \D,J)$.
\qed
\enddemo
\bigskip
The  theorem implies  the following  characterization
of  standard CR structures:
\proclaim {Corollary 5.2}
  A compact homogeneous CR manifold
is standard if and only if a maximal connected compact group
of automorphisms has   1-dimensional center.
\endproclaim
\bigskip
\bigskip
\subhead 6.  Equivalences of homogeneous compact CR manifolds
\endsubhead
\bigskip
The goal of this section is to determine
when two simply connected homogeneous compact CR manifolds
$(M=G/L, \D, J), \, (M'=G'/L', \D', J')$ are CR diffeomorphic. We will give our
results considering the cases of standard and non-standard CR manifolds
separately.\par
\bigskip
\subsubhead 6.1 The case of non-standard CR manifolds
\endsubsubhead
\medskip
Remark that a non-standard CR manifold can not be CR diffeomorphic to a standard one
(see, for instance,  Corollary 4.2). Moreover, two distinct
non-standard CR manifolds $G/L$ in Table 2, with $G/L \neq \text{Spin}_7/SU_3$,
are not CR diffeomorphic,  because in this case,
 $G$ coincides with a maximal connected
compact semisimple group $A^{ss}(M,\D,J)$ of automorphisms of $M = G/L$.
Moreover, by Theorem 4.5 (2),  for any invariant non-standard CR structure $(\D,J)$ on
 $M =\text{Spin}_7/SU_3$,  we have that
$A^{ss}(M,\D,J) = SO_8$ and that any other
 non-standard CR manifold $M' = G'/L' \neq \text{Spin}_7/SU_3$
 is CR diffeomorphic to $M$ if and if
$M' = SO_8/SO_6$. In particular, this implies that in order
to determine all compact CR homogeneous manifolds $(M' = G'/L', \D', J')$,
which are CR equivalent to a given compact homogeneous CR manifold $(M = G/L, \D, J)$,
there is no loss of generality if one assumes that $G \neq \text{Spin}_7$ and
that $G'/L' = G/L$.\par
Finally, it is known (see [AS]) that invariant CR structures $(\D,J)$ (considered
up to sign of $J$) on a given non-standard homogeneous
CR manifold $M = G/L$ with a fixed contact structure $\D$,
are naturally parameterized by points of the unite disc $D \subset \Bbb C$.
\par
It remains to find out   when two of such CR structures are CR equivalent. The
answer is given in the following proposition.\par
\bigskip
\proclaim{Proposition 6.1} Let $(M = A/B, \D)$ be a
 homogeneous contact manifold of a connected compact semisimple Lie group
$A \neq \text{Spin}_7$, which admits an $A$-invariant non-standard CR structure $(\D, J)$
(see Table 2).  Let also $(\D, J_t)$ and $(\D, J_{t'})$ be two
invariant CR structure on $M = A/B$ of the family
of non-standard CR structures parametrized by the points $t\in D\setminus\{0\} \subset \C$
of the punctured unit disc  in $\C$, as described in Cor. 5.2 and Prop. 6.3 and 6.4 in [AS].\par
Then $(\D, J_t)$ is CR equivalent to $(\D, J_{t'})$ (up to sign of
$J_{t'}$) if and only if
$|t| = |t'|$. \par
\endproclaim
\demo{Proof} Let
$$\phi:  (M = A/B, \D, J_t) \to (M = A/B, \D, J_{t'})$$
be a CR diffeomorphism. By Theorems 4.5 and 5.1, $A$ is a maximal compact
group of automorphisms in $Aut(M, \D, J_t)$ and in $Aut(M, \D, J_{t'})$.
Therefore $\phi$ transforms $A$ into $A' = \phi\circ A \circ \phi^{-1}$,
which is a maximal
compact group of automorphisms of $(M, \D, J_{t'})$. Since any two
maximal compact subgroups of  $Aut(M, \D, J_{t'})$ are conjugated,
without loss of generality we may assume that $A' = A$. Hence
$\phi$ induces a Lie group automorphism of $A$
which preserves the isotropy $B$, the contact element $Z$ associated with $\D$
(up to a scaling)
and transforms the holomorphic subspace $\m^{10}_t$ into
$\m^{10}_{t'}$. This means that $(\D, J_t)$ and $(\D, J_{t'})$, with $t,t' \in D \subset \C$,
are
CR equivalent if and only if  there exists a Lie automorphism
$\phi: \a^\C \to \a^\C$ such that
\roster
\item"i)" $\phi(\a) = \a$ and $\phi(\b) = \b$;
\item"ii)" $\phi(Z) \in \R Z$;
\item"iii)" $\phi(\m^{10}_t)  = \m^{10}_{t'}$.
\endroster
First  we consider an  inner automorphism $\phi$ which satisfies i), ii) and iii).
In this case, we may assume that it is of the form
$\phi =  \operatorname{Ad}_{\exp(s Z)}$. Indeed, $\phi$ verifies
 i) and ii) if and only if it is of the form
$\phi =  \operatorname{Ad}_{\exp(s Z)}\cdot \operatorname{Ad}_{\exp(X)}$,
for some $X\in \b$ and $s \in \R$. Since $\m^{10}_t$ is
$\operatorname{Ad}_\b$-invariant, we may neglect the factor $\operatorname{Ad}_{\exp(X)}$, which preserves
$\m^{10}_t$. \par
Using the explicit description of $Z$ and of $\m^{10}_t$
in Cor. 5.2 and Prop. 6.3 and 6.4 in \cite{AS}, the reader can
easily check that $\operatorname{Ad}_{\exp(sZ)}$ acts on
the space $\m^{10}_t$ by transforming it into the space $\m^{10}_{t'}$,
where $t' = e^{iCs} t$ for some constant $C\neq 0$ depending only on the
Lie algebra $\a$. This shows that there exists an inner automorphism which
verifies i), ii) and iii) if and only if $|t| = |t'|$.\par
It remains to consider the case when $\phi$ is an outer automorphism.
Composing it with an inner automorphism, we may always assume that
it preserves a Cartan subalgebra $\h = (\h\cap \b) + \R Z \subset \a$.\par
By Cor. 5.2, Prop. 6.3 or of Prop. 6.4 in \cite{AS}, we know that in all cases
of Table 2, there exists either one or two pairs of equivalent $\b$-moduli
in $\m^\C$. Assume for simplicity that there is only one pair $(\m_1,
\m_2)$ of equivalent $\b$-moduli. We may also assume that each $\m_i$, $i = 1,2$,
is also a $(\b + \R Z)$-module with highest weight $\alpha_i$, where
$\alpha_i$ is a root. Then
$\m^{10}_t$  and $\m^{10}_{t'}$ are two irreducible $\b$-moduli with highest weight vectors
$E_{\alpha_1} + t E_{\alpha_2}$ and $E_{\alpha_1} + t' E_{\alpha_2}$,
respectively,  where $E_{\alpha_i}$ is the root vector
with root $\alpha_i$ in the Chevalley normalization.\par
Since $\phi$ preserves the root system of $(\a, \h)$, either $\phi$
preserves the moduli $\m_i$ or interchange them. In the first case
$\phi(\m^{10}_t) = \m^{10}_t$ and in the second case
$\phi(\m^{10}_t) = \m^{10}_{1/t}$. Similar arguments show that the same conclusion
holds  also when
there exists two pairs of equivalent
$\b$-moduli in $\m^\C$.\par
Since $1/t \notin D$, it follows that an outer automorphism
verifies i), ii) and iii) if and only if $t = t'$. \qed
\enddemo
\bigskip
\subsubhead 6.2 The case of standard CR manifolds
\endsubsubhead
\medskip
 Now  we consider  the
standard homogeneous CR manifolds. Enlarging the group $G$ of automorphisms
of a standard homogeneous
CR manifold $(M=G/L, \D, J)$, we may always assume that $G = A^{ss}(M,\D,J)$.\par
Assume that  $(M=G/L, \D, J)$ and $(M'=G'/L', \D', J')$
are standard and CR diffeomorphic. Then, using a CR diffeomorphism, we may identify $M'$ with $M$ and
$G' = A^{ss}(M',\D',J')$ with a transitive subgroup of $\text{Aut}(M,\D,J)$ which is
conjugated to $G$. Using this conjugation, we may also identify $G$ with $G'$ and  $L$ with
$L'$. Therefore, the problem reduces to the description of all pairs of
standard invariant CR structure
$(\D,J) $ and $(\D',J')$ on the same homogeneous manifold $M= G/L$ that are CR equivalent.
\par
 The following proposition gives a necessary condition
for two  invariant standard CR structures on a
given homogeneous manifold $M = G/L$ to be CR equivalent.\par
\bigskip
\proclaim{Proposition 6.2} Assume that  $(\D,J)$ and $(\D',J')$
are two invariant standard  CR structures on a homogenous
CR manifold $M = G/L$. If they are
CR equivalent  then
the associated flag manifolds $F=G/K $ and
$F'=G/K'$ are biholomorphic
with respect to the invariant complex structures $J_F$ and  $J_{F'}$, induced by $J$ and $J'$,
respectively.
\endproclaim
\demo{Proof} By Corollary 3.3,
we may assume that $G$ is equal to a maximal semisimple groups of
CR transformations of $M$. Therefore, if we denote by $Z$ and $Z'$ two contact
elements associated with $\D$ and $\D'$, then
$(\D, J)$ and $(\D', J')$ are CR equivalent
only if there is a Lie  automorphism $\phi: \g = Lie(G) \to \g = Lie(G)$
 with the following properties:
$$\text{a)}\ \ \phi(\l) = \l\ ;\qquad
\text{b)}\ \ \phi(\R Z) = \R Z'\ ;\qquad
 \text{c)}\ \ \phi(\m^{10}) = \m^{10}{}'\ .$$
Since a) and b) imply that $\phi(\k) = \k'$, the automorphism $\phi$ induces
a $G$-equivariant biholomorphic map between $F$ and $F'$, with respect to the complex structures
associated with $\m^{10}$ and $\m^{10}{}'$.
\qed
\enddemo
Now we describe all  invariant CR structures $(M=G/L,\D, J)$ with given
associated flag manifold $(F=G/K,J_F)$ up to  CR equivalence.
 Let $\h$ be a Cartan subalgebra of $\k^{\Bbb C}$ and
$R_K$ and $R$  the root systems of $\k^{\Bbb C}$ and of
$\g^{\Bbb C}$  with respect to $\h$. Fix a basis $\Pi_K$ of $R_K$.
It is known that there exists a 1-1 correspondence
between
invariant complex structures $J_F$ on $F=G/K$ and   bases $\Pi$ of $R$
containing  $\Pi_K$. We call such a root system $R$ and such a basis
$\Pi$ {\it a root system and a basis adapted to the flag manifold\/} \par
Recall that any adapted  basis
can be represented by  Dynkin graph with black and white nodes,
where the black nodes are associated with the simple roots
in $\Pi' = \Pi \setminus \Pi_K$. \par
The complex structure which is associated with an adapted basis $\Pi$ of the above kind is determined
by the holomorphic subspace
$$ \m^{10} = \g(\Pi\setminus \Pi_K)\ ,$$
where we use the notation $\g(S)$, $S \subset R$, to denote the  subalgebra generated
by the root vectors $E_\alpha$, $\alpha \in S$.\par
\bigskip
Now, let us fix a flag manifold $(F = G/K, J_F)$ with an
invariant complex structure $J_F$. If $\Pi$ is the corresponding
adapted basis, we will
denote  by $\pi_1, \dots , \pi_m$ the fundamental weights, which corresponds
to the `black' simple
roots $\Pi' = \Pi \setminus \Pi_K = \{ \alpha_1, \dots , \alpha_m \}$.
We call $\pi_1, \dots \pi_m$ the  {\it black weights associated with the adapted
basis $\Pi$\/}.
\par
\medskip
\proclaim{Theorem 6.3} Let $(F=G/K,J_F)$ be a flag manifold and
let $\pi_1, \dots \pi_m$ the  black weights associated with the
adapted basis $\Pi$ corresponding to $J_F$. Then there exists a
1-1 correspondence between standard homogeneous CR manifolds $(M
=G/L,\D,J)$  up to  CR equivalence and the set of all m-tuples
$\bar p = (p_1, \dots , p_m )\in \Bbb Z^m $ such that the $p_i
\neq 0 , i = 1,\dots , m $  have  no common divisor.\par
 Any such
m-tuple $\bar p$ corresponds to the  following homogeneous CR
manifold $(M=G/L, \D_Z, J)$: \roster
\item
 the subgroup
 $L$ is the connected subgroup of $G$, generated by
the subalgebra
$$\l = [\k,\k] + Z(\k) \cap \text{ker}\ \theta\ ,\quad  \text{where}\ \  \theta = p_1 \pi_1 + \dots + p_m \pi_m\ ;$$
\item the contact structure $\D $ is defined by the contact element $Z = \B^{-1}\theta$;
\item the
CR structure $J$ is the one associated with the holomorphic subspace $\m^{10} = \g(\Pi')$,
where $\Pi' = \{ \alpha_1, \dots, \alpha_m\}$ are the
`black'  roots of the adapted basis $\Pi$.
\item The Levi form $d \theta \circ J$ is positively defined if
and only if all $p_i$ are positive.
\endroster
\endproclaim
\demo{Proof} It follows immediately from Proposition 6.2 and Proposition 2.1.\qed
\enddemo

\enddocument

\bye